\documentclass{llncs}
\usepackage[utf8]{inputenc}

\RequirePackage{amsmath}
\usepackage{amsmath, bm, amsfonts}
\usepackage{graphicx}
\usepackage{wrapfig}

\usepackage{xcolor}

\usepackage{hyperref} 
\hypersetup{
colorlinks=false, 
breaklinks=true,
urlcolor= blue, 
linkcolor= blue,
citecolor=red, 
}

\newcommand{\ipujson}[0]{{\sc bqpjson}}

\newcommand{\ipig}[0]{{\sc dwig}}

\newcommand{\ipusolvers}[0]{{\sc bqpsolvers}}

\newcommand{\timeout}[0]{TO}

\begin{document}

\title{Evaluating Ising Processing Units\\ with Integer Programming}

%\author{}
\author{Carleton Coffrin \and Harsha Nagarajan \and Russell Bent}
\institute{Los Alamos National Laboratory, Los Alamos, New Mexico, USA}

\date{}
%\date{Received: date / Accepted: date}
% The correct dates will be entered by the editor

\maketitle

\begin{abstract}
The recent emergence of novel computational devices, such as adiabatic quantum computers, CMOS annealers, and optical parametric oscillators, present new opportunities for hybrid-optimization algorithms that are hardware accelerated by these devices. In this work, we propose the idea of an Ising processing unit as a computational abstraction for reasoning about these emerging devices. The challenges involved in using and benchmarking these devices are presented and commercial mixed integer programming solvers are proposed as a valuable tool for the validation of these disparate hardware platforms.  The proposed validation methodology is demonstrated on a D-Wave 2X adiabatic quantum computer, one example of an Ising processing unit. The computational results demonstrate that the D-Wave hardware consistently produces high-quality solutions and suggests that as IPU technology matures it could become a valuable co-processor in hybrid-optimization algorithms.
\end{abstract}

\section{Introduction}
\label{sec:intro}

As the challenge of scaling traditional transistor-based Central Processing Unit (CPU) technology continues to increase, experimental physicists and high-tech companies have begun to explore radically different computational technologies, such as adiabatic quantum computers (AQCs) \cite{Johnson2011}, gate-based quantum computers \cite{ibm_quantum,45919,chmielewski2018cloud}, CMOS annealers \cite{7063111,6662276,fuhitsu_da}, neuromorphic computers \cite{ibm_nuro,8259423,schuman2017survey}, memristive circuits \cite{1712.07046,1808.09999}, and optical parametric oscillators \cite{mcmahon2016fully,Inagaki603,7738704}.  The goal of all of these technologies is to leverage the dynamical evolution of a physical system to perform a computation that is challenging to emulate using traditional CPU technology (e.g., the simulation of quantum physics) \cite{Feynman1982-FEYSPW}.  Despite their entirely disparate physical implementations, AQCs, CMOS annealers, memristive circuits, and optical parametric oscillators are unified by a common mathematical abstraction known as the Ising model, which has been widely adopted by the physics community for the study of naturally occurring discrete optimization processes \cite{RevModPhys.39.883}.  Furthermore, this kind of ``Ising machine'' \cite{mcmahon2016fully,Inagaki603} is already commercially available with more than 2000 decision variables in the form of AQCs developed by D-Wave Systems \cite{dwave_customers}.

The emergence of physical devices that can quickly solve Ising models is particularly relevant to the constraint programming, artificial intelligence and operations research communities, because the impetus for building these devices is to perform discrete optimization.  As this technology matures, it may be possible for this specialized hardware to rapidly solve challenging combinatorial problems, such as Max-Cut \cite{Haribara2016} or Max-Clique \cite{10.3389/fphy.2014.00005}.  Preliminary studies have suggested that some classes of Constraint Satisfaction Problems may be effectively encoded in such devices because of their combinatorial structure \cite{10.3389/fict.2016.00014,10.3389/fphy.2014.00056,Rieffel2015,1506.08479}.
Furthermore, an Ising model coprocessor could have significant impacts on solution methods for a variety of fundamental combinatorial problem classes, such as MAX-SAT \cite{deGivry2003,Morgado2013,McGeoch:2013:EEA:2482767.2482797} and integer programming \cite{Nieuwenhuis2014}. At this time, however, it remains unclear how established optimization algorithms should leverage this emerging technology.  This paper helps to address this gap by highlighting the key concepts and hardware limitations that an algorithm designer needs to understand to engage in this emerging and exciting computational paradigm.

Similar to an arithmetic logic unit (ALU) or a graphics processing unit (GPU), this work proposes the idea of an Ising processing unit (IPU) as the computational abstraction for wide variety of physical devices that perform optimization of Ising models.  This work begins with a brief introduction to the IPU abstraction and its mathematical foundations in Section \ref{sec:ising}.  Then the additional challenges of working with real-world hardware are discussed in Section \ref{sec:ipu} and an overview of previous benchmarking studies and solution methods are presented in Section \ref{sec:benchmarking}.  Finally, a detailed benchmarking study of a D-Wave 2X IPU is conducted in Section \ref{sec:computations}, which highlights the current capabilities of such a device.  The contributions of this work are as follows,
\begin{enumerate}
    \item The first clear and concise introduction to the key concepts of Ising models and the limitations of real-world IPU hardware, both of which are necessary for optimization algorithm designers to effectively leverage these emerging hardware platforms (Section \ref{sec:ising} and Section \ref{sec:ipu}).
    \item Highlighting that integer programming has been overlooked by recent IPU benchmarking studies (Section \ref{sec:benchmarking}), and demonstrating the value of integer programming for filtering easy test cases (Section \ref{sec:computations:cases}) and verifying the quality of an IPU on challenging test cases (Section \ref{sec:computations:ipu}).
\end{enumerate}
%
%
% (1) presenting foundational Ising model results from other disciplines in terminology that is familiar to the optimization community (Section \ref{sec:ising}), and (2) highlighting key features of real-world IPUs, to provide context for future algorithmic developments utilizing these analog devices (Section \ref{sec:ipu}).  The technical contributions include (1) developing a benchmarking methodology for IPUs, which is enabled by proposed open-source software tools (Section \ref{sec:tools});  (2) demonstrating the proposed benchmarking tools by conducting a baseline evaluation of a D-Wave 2X IPU (Section \ref{sec:computations}); and (3) illustrating that 
%
% To the best of our knowledge, this is the first benchmarking study to use an integer programming solver for identifying challenging IPU test cases.
%
Note that, due to the maturity and commercial availability of the D-Wave IPU, this work often refers to that architecture as an illustrative example.  However, the methods and tools proposed herein are applicable to all emerging IPU hardware realizations, to the best of our knowledge.

% Old Paragraph.
% To assist optimization researchers in exploring how to integrate this novel hardware into established algorithms, in this work, we propose the idea of Ising Processing Units (IPUs) as a computational abstraction for any physical device that performs optimization of an Ising Model.  We review a number of unique features of working with current IPU devices and discuss a number of unexpected challenges around benchmarking this kind of hardware (Section \ref{sec:ipu}).  We then propose a number of software tools to help address some of the benchmarking challenges (Section \ref{sec:tools}) and perform a baseline study of a D-Wave adiabatic quantum computer using the proposed benchmarking methodology (Section \ref{sec:computations}).  Due to the maturity and commercial availability of the D-Wave IPU, this work often refers to that architecture as an illustrative example.  However, based on our understanding of other IPU technologies, we anticipate that the ideas proposed herein are applicable to those architectures as well.

\section{A Brief Introduction to Ising Models}
\label{sec:ising}

This section introduces the notations of the paper and provides a brief introduction to Ising models, the core mathematical abstraction of IPUs.
The Ising model refers to the class of graphical models where the nodes, ${\cal N}$, represent {\em spin} variables (i.e., $\sigma_i \in \{-1,1\} ~\forall i \in {\cal N}$) and the edges, ${\cal E}$, represent {\em interactions} of spin variables (i.e., $\sigma_i \sigma_j ~\forall i,j \in {\cal E}$).  A local {\em field} $\bm h_i ~\forall i \in {\cal N}$ is specified for each node, and an interaction strength $\bm J_{ij} ~\forall i,j \in {\cal E}$ is specified for each edge.  Given these data, the {\em energy} of the Ising model is defined as,
\begin{align}
    E(\sigma) &= \sum_{i,j \in {\cal E}} \bm J_{ij} \sigma_i \sigma_j + \sum_{i \in {\cal N}} \bm h_i \sigma_i \label{eq:ising_eng}
\end{align}
Applications of the Ising model typically consider one of two tasks. First, some applications focus on finding the lowest possible energy of the Ising model, known as a {\em ground state}.  That is, finding the globally optimal solution of the following binary quadratic optimization problem:
\begin{align}
    & \min: E(\sigma) \label{eq:ising_opt}\\
    & \mbox{s.t.: } \sigma_i \in \{-1, 1\}  ~\forall i \in {\cal N} \nonumber
\end{align}
Second, other applications are interested in sampling from the Boltzmann distribution of the Ising model's states:
\begin{align}
    & Pr(\sigma) \propto \bm e^{\frac{-E(\sigma)}{\bm \tau}} \label{eq:ising_smpl}
\end{align}
where $\bm \tau$ is a parameter representing the {\em effective temperature} of the Boltzmann distribution \cite{Zdeborova2016statistical}.
It is valuable to observe that in the Boltzmann distribution, the lowest energy states have the highest probability.  Therefore, the task of sampling from a Boltzmann distribution is similar to the task of finding the lowest energy of the Ising model.  Indeed, as $\bm \tau$ approaches 0, the sampling task smoothly transforms into the aforementioned optimization task. 
This paper focuses exclusively on the mathematical program presented in \eqref{eq:ising_opt}, the optimization task.

\paragraph{Frustration:}
The notion of frustration is common in the study of Ising models and refers to any instance of \eqref{eq:ising_opt} where the optimal solution, $\sigma^*$, satisfies the property,
%\begin{subequations}
\begin{align}
    E(\sigma^*) > \sum_{i,j \in {\cal E}} - |\bm J_{ij}| - \sum_{i \in {\cal N}} |\bm h_i| 
\end{align}
%\end{subequations}
A canonical example is the following three node problem: 
%
%\begin{subequations}
\begin{align}
    \bm h_{1} = 0, ~\bm h_{2} &= 0, ~\bm h_{3} = 0, ~\bm J_{12} = -1, ~\bm J_{23} = -1, ~\bm J_{13} = 1
\end{align}
%\end{subequations}
%
Observe that, in this case, there are a number of optimal solutions such that $E(\bm \sigma^*) = -2$ but none such that $E(\sigma) = \sum_{i,j \in {\cal E}} - |\bm J_{ij}| = -3$.  Note that frustration has important algorithmic implications as greedy algorithms are sufficient for optimizing Ising models without frustration.

\paragraph{Gauge Transformations:}
A valuable property of the Ising model is the gauge transformation, which characterizes an equivalence class of Ising models.  For illustration, consider the optimal solution of Ising model $S$, $\bm \sigma^{s*}$.  One can construct a new Ising model $T$ where the optimal solution is the same, except that $\bm \sigma^{t*}_i = - \bm  \sigma^{s*}_i$ for a particular node $i \in {\cal N}$ is as follows:
\begin{subequations}
\begin{align}
    \bm J^t_{ij} &= - \bm J^s_{ij}  ~\forall i,j \in {\cal E}(i) \\
    \bm h^t_i &= - \bm h^s_i
\end{align}
\end{subequations}
where ${\cal E}(i)$ indicates the neighboring edges of node $i$.  This $S$-to-$T$ manipulation is referred to as a gauge transformation.  Given a complete source state $\bm \sigma^{s}$ and a complete target state $\bm \sigma^{t}$, this transformation is generalized to all of $\sigma$ by,
% s^s_i s^s_j |  s^t_i s^t_j | J sign
% -1 -1 | -1 -1 same
% -1 -1 |  1 -1 inv.
% -1 -1 | -1  1 inv.
% -1 -1 |  1  1 same
%
%  1 -1 | -1 -1 inv.
%  1 -1 |  1 -1 same
%  1 -1 | -1  1 same
%  1 -1 |  1  1 inv.
%
% -1  1 | -1  1 same
%...
%
%  1  1 |  1  1 same
%...
%
\begin{subequations}
\begin{align}
    \bm J^t_{ij} &= \bm J^s_{ij} \bm \sigma^s_i \bm \sigma^s_j \bm \sigma^t_i \bm \sigma^t_j ~\forall i,j \in {\cal E} \\
    \bm h^t_i &= \bm h^s_i \bm \sigma^s_i \bm \sigma^t_i ~\forall i \in {\cal N}
\end{align}
\end{subequations}
It is valuable to observe that by using this gauge transformation property, one can consider the class of Ising models where the optimal solution is $\sigma^{*}_i = -1 ~\forall i \in {\cal N}$ or any arbitrary vector of ${-1,1}$ values without loss of generality.

\paragraph{Bijection of Ising and Boolean Optimization:}
It is also useful to observe that there is a bijection between Ising optimization (i.e., $\sigma \in \{-1,1\}$) and Boolean optimization (i.e., $x \in \{0,1\}$).  The transformation of $\sigma$-to-$x$ is given by,
\begin{subequations} \label{eq:spin2bool}
\begin{align}
    \sigma_i &= 2x_i - 1  ~\forall i \in {\cal N} \\ 
    \sigma_i\sigma_j &= 4x_ix_j - 2x_i - 2x_j + 1  ~\forall i,j \in {\cal E}
\end{align} \label{}
\end{subequations} 
and the inverse $x$-to-$\sigma$ is given by,
\begin{subequations}
\begin{align}
    x_i &= \frac{\sigma_i + 1}{2}  ~\forall i \in {\cal N} \\ 
    x_i x_j &= \frac{\sigma_i \sigma_j + \sigma_i + \sigma_j + 1}{4}  ~\forall i,j \in {\cal E} %\\
\end{align}
\end{subequations}
Consequently, any results from solving Ising models are also immediately applicable to the following class of Boolean optimization problems:
%
%\begin{subequations}
\begin{align}
    & \min: \sum_{i,j \in {\cal E}} \bm c_{ij} x_i x_j + \sum_{i \in {\cal N}} \bm c_i x_i \label{eq:boolean_opt} \\
    & \mbox{s.t.: } x_i \in \{0, 1\}  ~\forall i \in {\cal N} \nonumber
\end{align}
%\end{subequations}
%
The Ising model provides a clean mathematical abstraction for understanding the computation that IPUs perform.  However, in practice, a number of hardware implementation factors present additional challenges for computing with IPUs.

\section{Features of Analog Ising Processing Units}
\label{sec:ipu}

The core inspiration for developing IPUs is to take advantage of the natural evolution of a discrete physical system to find high-quality solutions to an Ising model \cite{Johnson2011,mcmahon2016fully,6662276,1712.07046}.  Consequently, to the best of our knowledge, all IPUs developed to date are analog machines, which present a number of challenges that the optimization community is not accustomed to considering.

\paragraph{Effective Temperature:}
The ultimate goal of IPUs is to solve the optimization problem \eqref{eq:ising_opt} and determine the globally optimal solution to the input Ising model.  In practice, however, a variety of analog factors preclude IPUs from reliably finding globally optimal solutions.  As a first-order approximation, current IPUs behave like a Boltzmann sampler \eqref{eq:ising_smpl} with some hardware-specific effective temperature, $\bm \tau$ \cite{dwave_boltzmann}.  It has also been observed that the effective temperature of an IPU can vary around a nominal value based on the Ising model that is being executed \cite{PhysRevA.94.022308}.  This suggests that the IPU's performance can change based on the structure of the problem input.

\paragraph{Environmental Noise:}
One of the primary contributors to the sampling nature of IPUs are the environmental factors.  All analog machines are subject to faults due to environmental noise; for example, even classical computers can be affected by cosmic rays.  However, given the relative novelty of IPUs, the effects of environmental noise are noticeable in current hardware.  The effects of environmental noise contribute to the perceived effective temperature $\bm \tau$ of the IPU.

\paragraph{Coefficient Biases:}
Once an Ising model is input into an IPU, its coefficients are subject to at least two sources of bias.  The first source of bias is a model programming error that occurs independently each time the IPU is configured for a computation.  This bias is often mitigated by programming the IPU multiple times with an identical input and combining the results from all executions.  The second source of bias is a persistent coefficient error, which is an artifact of the IPU manufacturing and calibration process.
 Because this bias is consistent across multiple IPU executions, this source of bias is often mitigated by performing multiple gauge transformations on the input and combining the results from all executions.

\begin{figure}[t]
    \begin{center}
    \includegraphics[scale=0.90]{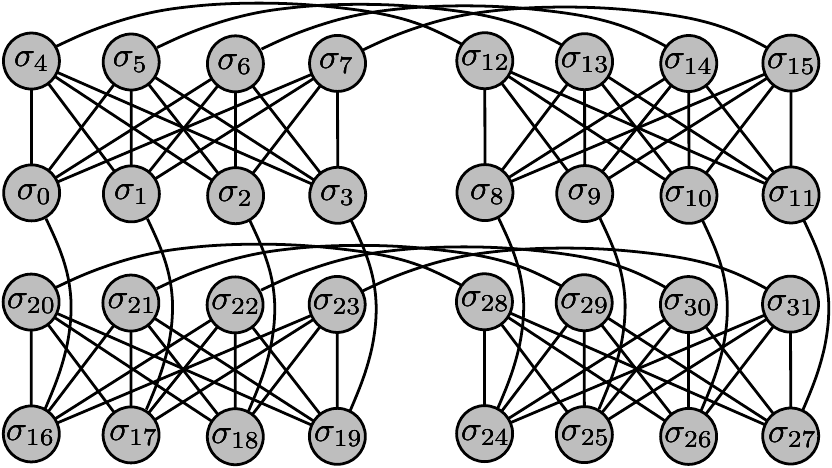}
    \end{center}
    \vspace{-0.4cm}
    \caption{A 2-by-2 Chimera Graph Illustrating the Variable Product Limitations of a D-Wave 2X IPU.}
    \vspace{-0.2cm}
    \label{fig:chimera}
\end{figure}

\paragraph{Problem Coefficients:}
In traditional optimization applications, the problem coefficients are often rescaled to best suit floating-point arithmetic.  Similarly, IPUs have digital-to-analog converters that can encode a limited number of values; typically these values are represented as numbers in the range of -1 to 1.  Some IPUs allow for hundreds of steps within this range, \cite{Johnson2011,6662276} whereas others support only the discrete set of \{-1, 0, 1\} \cite{mcmahon2016fully}.  In either case, the mathematical Ising model must be rescaled into the IPU's operating range.  However, this mathematically equivalent transformation can result in unexpected side effects because the coefficients used in the IPU hardware are perturbed by a constant amount of environmental noise and hardware bias, which can outweigh small rescaled coefficient values.

\paragraph{Topological Limitations:}
Another significant feature of IPUs is a restricted set of variable products.  In classical optimization  (e.g., \eqref{eq:ising_opt}), it is assumed that every variable can interact with every other variable, that is, an Ising model where an edge connects every pair of variables.  However, because of the hardware implementation of an IPU, it may not be possible for some variables to interact.  For example, the current D-Wave IPUs are restricted to the {\em chimera} topology, which is a two-dimensional lattice of {\em unit cells}, each of which consist of a 4-by-4 bipartite graph (e.g., see Figure \ref{fig:chimera}).  In addition to these restrictions, fabrication errors can also lead to random failures of nodes and edges in the IPU hardware.  Indeed, as a result of these minor imperfections, every D-Wave IPU developed to date has a unique topology \cite{Boixo2014,PhysRevX.6.031015,king2015benchmarking}.  Research and development of algorithms for embedding various kinds of Ising models into a specific IPU topology is still an active area of research \cite{10.3389/fict.2016.00014,Boothby:2016:FCM:2877060.2877142,1406.2741,Klymko2014}.

\subsection{Challenges of Benchmarking Ising Processing Units}

These analog hardware features present unique challenges for benchmarking IPUs that fall roughly into three categories: (1) comparing to established benchmark libraries; (2) developing Ising model instance generators for testing and; (3) comparing with classical optimization methods.

\paragraph{Benchmark Libraries:}
Research and development in optimization algorithms has benefited greatly from standardized benchmark libraries \cite{miplib2010,Gent1999,Hoos00satlib:an}.  However, direct application of these libraries to IPUs is out of scope in the near term for the following reasons: (1) the Ising model is a binary quadratic program, which is sufficiently restrictive to preclude the use of many standard problem libraries; (2) even in cases where the problems of interest can be mapped directly to the Ising model (e.g., Max-Cut, Max-Clique), the task of embedding given problems onto the IPU's hardware graph can be prohibitive \cite{coffrin2016challenges}; and (3) even if an embedding can be found, it is not obvious that the problem's coefficients will be amenable to the IPU's operating range.

\paragraph{Instance Generation Algorithms:}
%
%Another option is to build a new standardized benchmark library specifically for the evaluation of IPUs.  However, given that every IPU is unique and within an IPU each variable also has individual properties (e.g., persistent bias), it is not immediately clear that a new static benchmark library could be applied to multiple IPUs.
%
Due to these challenges, the standard practice in the literature is to generate a collection of instances for a given IPU and use these cases for the evaluation of that IPU \cite{king2015benchmarking,PhysRevA.92.042325,1701.04579,PhysRevX.6.031015}. The hope being that these instances provide a reasonable proxy for how real-world applications might perform on such a device. 

%The hope being that if a sufficient number of test instances are generated using the same algorithm, the results will be replicable across multiple IPUs.

\paragraph{Comparison with Classical Algorithms:}
Because of the radically different hardware of CPUs vs IPUs and the stochastic nature of the IPUs, conducting a fair comparison of these two technologies is not immediately clear \cite{PhysRevA.94.022337,1703.00622,1701.04579}.  
%For example, is the comparison more fair or less fair if classical algorithms utilize multiple cores or GPUs?  
Indeed, comparisons of D-Wave's IPU with classical algorithms have resulted in vigorous discussions about what algorithms and metrics should be used to make such comparisons \cite{aaronson_blog1,king2015benchmarking,aaronson_blog2}.  It is widely accepted that IPUs do not provide optimality guarantees and are best compared to heuristic methods (e.g. local search) in terms of runtime performance.  
%However, it is less clear whether heuristics should be specialized to solve problems for a specific IPU architecture (e.g., \cite{1409.3934,HFS_impl_2017}) or whether classical methods should be designed for the  most general form of the Ising model \eqref{eq:ising_opt}.
This debate will most likely continue for several years.  In this work, our goal is not to answer these challenging questions but rather to highlight that commercial mixed integer programming solvers are valuable and important tools for exploring these questions.

\section{A Review of Ising Processing Unit Benchmarking Studies}
\label{sec:benchmarking}

Due to the challenges associated with mapping established optimization test cases to specific IPU hardware \cite{coffrin2016challenges}, the IPU benchmarking community has adopted the practice of generating Ising model instances on a case-by-case basis for specific IPUs \cite{king2015benchmarking,PhysRevA.92.042325,1701.04579,PhysRevX.6.031015} and evaluating these instances on a variety of solution methods.  The following subsections provide a brief overview of the instance generation algorithms and solution methods that have been used in various IPU benchmarking studies.  The goals of this review are to: (1) reveal the lack of consistency across current benchmarking studies; (2) highlight the omission of integer programming methods in all of the recent publications and; (3) motivate the numerical study conducted in this work. 

\subsection{Instance Generation Algorithms}
\label{sec:tools_gen}

The task of IPU instance generation amounts to finding interesting values for $\bf h$ and $\bf J$ in \eqref{eq:ising_eng}.  In some cases the procedures for generating these values are elaborate \cite{PhysRevX.6.031015,king2015performance} and are designed to leverage theoretical results about Ising models \cite{PhysRevA.92.042325}.  A brief survey reveals five primary problem classes in the literature, each of which is briefly introduced.  For a detailed description, please refer to the source publication of the problem class.

\paragraph{Random (RAN-k and RANF-k):}
To the best of our knowledge, this general class of problem was first proposed in \cite{McGeoch:2013:EEA:2482767.2482797} and was later refined into the RAN-k problem in \cite{king2015benchmarking}.  The RAN-k problem consists simply of assigning each value of $\bf h$ to 0 and each value of $\bf J$ uniformly at random from the set
\begin{align}
\{ -\bm k, -\bm k+1, \dots, -2, -1, 1, 2, \dots, \bm k-1, \bm k \} \label{eq:k_set}
\end{align}
The RANF-k problem is a simple variant of RAN-k where the values of $\bf h$ are also selected uniformly at random from \eqref{eq:k_set}.  As we will later see, RAN-1 and RANF-1, where $\bm h, \bm J \in \{-1, 1\}$, are an interesting subclass of this problem.

\paragraph{Frustrated Loops (FL-k and FCL-k):}
The frustrated loop problem was originally proposed in \cite{PhysRevA.92.042325} and then later refined to the FL-k problem in \cite{king2015performance}.  It consists of generating a collection of random cycles in the IPU graph.  In each cycle, all of the edges are set to $-1$ except one random edge, which is set to $1$ to produce {\em frustration}.  A scaling factor, $\alpha$, is used to control how many random cycles should be generated, and the parameter $k$ determines how many cycles each edge can participate in.  A key property of the FL-k generation procedure is that two globally optimal solutions are maintained at $\sigma_i = -1 ~\forall i \in {\cal N}$ and $\sigma_i = 1 ~\forall i \in {\cal N}$ \cite{king2015performance}.  However, to obfuscate this solution, a gauge transformation is often applied to make the optimal solution a random assignment of $\sigma$.

A variant of the frustrated loop problem is the frustrated {\em cluster} loop problem, FCL-k \cite{1701.04579}.  The FCL-k problem is inspired by the chimera network topology (i.e., Figure \ref{fig:chimera}).  The core idea is that tightly coupled variables (e.g., $\sigma_0...\sigma_7$ in Figure \ref{fig:chimera}) should form a {\em cluster} where all of the variables take the same value.  This is achieved by setting all of the values of $\bm J$ within the cluster to $-1$.  For the remaining edges between clusters, the previously described frustrated cycles generation scheme is used.  Note that a polynomial time algorithm is known for solving the FCL-k problem class on chimera graphs \cite{1703.00622}.

It is worthwhile to mention that the FL-k and FCL-k instance generators are solving a cycle packing problem on the IPU graph.  Hence, the randomized algorithms proposed in \cite{PhysRevA.92.042325,1701.04579} are not guaranteed to find a solution if one exists. In practice, this algorithm fails for the highly constrained settings of $\alpha$ and $k$.

\paragraph{Weak-Strong Cluster Networks (WSCNs):}
The WSCN problem was proposed in \cite{PhysRevX.6.031015} and is highly specialized to the chimera network topology.  The basic building block of a WSCN is a pair of spin clusters in the chimera graph (e.g., $\sigma_0...\sigma_7$ and $\sigma_8...\sigma_{15}$ in Figure \ref{fig:chimera}).  In the {\em strong} cluster the values of $\bm h$ are set to the strong force parameter {\em sf} and in the {\em weak} cluster the values of $\bm h$ are set to the weak force parameter {\em wf}.  All of the values of $\bm J$ within and between this cluster pair are set to $-1$.  Once a number of weak-strong cluster pairs have been placed, the strong clusters are connected to each other using random values of $\bm J \in \{ -1, 1 \}$.  The values of {\em sf} = $-1.0$ and {\em wf} = 0.44 are recommended by \cite{PhysRevX.6.031015}.  The motivation for the WSCN design is that the clusters create deep local minima that are difficult for local search methods to escape.

\subsection{Solution Methods}

Once a collection of Ising model instances have been generated, the next step in a typical benchmarking study is to evaluate those instances on a variety of solution methods, including the IPU, and compare the results.  A brief survey reveals five primary solution methods in the literature, each of which is briefly introduced.  For a detailed description, please refer to the source publications of the solution method.

\paragraph{Simulated Annealing:}
The most popular staw-man solution method for comparison is Simulated Annealing \cite{Kirkpatrick671}.  Typically the implementation only considers a neighborhood of single variable flips and the focus of these implementations is on computational performance (e.g. using GPUs for acceleration).  The search is run until a specified time limit is reached.

\paragraph{Large Neighborhood Search:}
The state-of-the-art meta-heuristic for solving Ising models on the chimera graphs is a Large Neighborhood Search (LNS) method called the Hamze-Freitas-Selby (HFS) algorithm \cite{Hamze:2004:FT:1036843.1036873,1409.3934}.  The core idea of this algorithm is to extract low treewidth subgraphs of the given Ising model and then use dynamic programming to compute the optimal configuration of these subgraphs.  This extract and optimize process is repeated until a specified time limit is reached.  A key to this method's success is the availability of a highly optimized open-source C implementation \cite{HFS_impl_2017}.

\paragraph{Integer Programming:}
Previous works first considered integer quadratic programming \cite{McGeoch:2013:EEA:2482767.2482797} and quickly moved to integer linear programming  \cite{ibm_blog,1306.1202} as a solution method.  The mathematical programming survey \cite{Billionnet2007} provides a useful overview of the advantages and dis-advantages of various integer programming (IP) formulations.

Based on some preliminary experiments with different formulations, this work focuses on the following integer linear programming formulation of the Ising model, transformed into the Boolean variable space:
\begin{subequations}
\begin{align}
    & \min:  \sum_{i,j \in {\cal E}} \bm c_{ij} x_{ij} + \sum_{i \in {\cal N}} \bm c_i x_i + \bm c \\
    & \mbox{s.t.: } \nonumber \\
    & x_{ij} \geq x_{i} + x_{j}- 1, ~x_{ij} \leq x_{i}, ~x_{ij} \leq x_{j} ~\forall i,j \in {\cal E} \\
    %& 2 x_{ij} \geq x_{i} + x_{j}  ~\forall i,j \in {\cal E} \\
    %& x_{ij} \leq x_{i}, ~x_{ij} \leq x_{j} ~\forall i,j \in {\cal E} \\
    %& x_{ij} \leq x_{j}  ~\forall i,j \in {\cal E} \\
    & x_i \in \{0, 1\} ~\forall i \in {\cal N}, ~x_{ij} \in \{0, 1\} ~\forall i,j \in {\cal E} \nonumber 
    %& x_i \in \{0, 1\} ~\forall i \in {\cal N} \nonumber 
\end{align}
\end{subequations}
where the application of \eqref{eq:spin2bool} leads to,
\begin{subequations}
\begin{align}
    \bm c_{ij} &= \sum_{i,j \in {\cal E}} 4 \bm J_{ij} ~\forall i,j \in {\cal E} \\
    \bm c_{i} &= \sum_{i,j \in {\cal E}(i)} - 2 \bm J_{ij} + \sum_{i \in {\cal N}} 2 \bm h_i  ~\forall i \in {\cal N} \\
    \bm c &= \sum_{i,j \in {\cal E}} \bm J_{ij} - \sum_{i \in {\cal N}} \bm h_i 
\end{align}
\end{subequations}
In this formulation, the binary quadratic program defined in \eqref{eq:boolean_opt} is converted to a binary linear program by lifting the variable products $x_ix_j$ into a new variable $x_{ij}$ and adding linear constraints to capture the $x_{ij} = x_i \wedge x_j ~\forall i,j \in {\cal E}$ conjunction constraints.  Preliminary experiments of this work confirmed the findings of \cite{Billionnet2007}, that this binary linear program formulation is best on sparse graphs, such as the hardware graphs of current IPUs.

\paragraph{Adiabatic Quantum Computation:}
%
%\cite{Johnson2011} {\color{red} This reference seemse out of place?}
%
An adiabatic quantum computation (AQC) \cite{quant-ph-0001106} is a method for solving an Ising model via a quantum annealing process \cite{PhysRevE.58.5355}.  This solution method has two notable traits: (1) the AQC dynamical process features quantum tunneling \cite{Farhi472}, which can help it to escape from local minima; (2) it can be implemented in hardware (e.g. the D-Wave IPU).

%In this work, and related works, this solution method  assumes that a given Ising model will map directly on the hardware graph.  No attempt is made to transform a given Ising model into the specified IPU hardware (i.e. no embedding computation is necessary).  

\paragraph{Quantum Monte Carlo:}
Quantum Monte Carlo (QMC) is a probabilistic algorithm that can be used for simulating large quantum systems.  QMC is a very computationally intensive method \cite{9780792355526,PhysRevX.6.031015} and thus the primary use of QMC is not to compare runtime performance but rather to quantify the possible value of an adiabatic quantum computation that could be implemented in hardware at some point in the future.

\subsection{Overview}

To briefly summarize a variety of benchmarking studies, Table \ref{tbl:bench_overview} provides an overview of the problems and solution methods previous works have considered.  Although there was some initial interest in integer programming models \cite{McGeoch:2013:EEA:2482767.2482797,ibm_blog,1306.1202}, more recent IPU benchmark studies have not considered these solution methods and have focused exclusively on heuristic methods.  Furthermore, there are notable inconsistencies in the type of problems being considered.  As indicated by the last row in Table \ref{tbl:bench_overview}, the goal of this work is revisit the use of IP methods for benchmarking IPUs and to conduct a thorough and side-by-side study of all problem classes and solution methods proposed in the literature.  Note that, because this paper focuses exclusively on the quality and runtime of the Ising model optimization task \eqref{eq:ising_opt}, the study of SA and QMC are omitted as they provide no additional insights over the LNS \cite{king2015performance} and AQC \cite{PhysRevX.6.031015} methods. 

\begin{table*}[t]
\center
\begin{tabular}{|c||c|c|c|c|c||c|c|c|c|c|c|c|c|}
\hline
            & \multicolumn{5}{|c||}{Problem Classes} & \multicolumn{5}{|c|}{Solution Methods} \\
Publication & RAN & RANF & FL & FCL & WSCN & IP & SA & LNS & QMC & AQC \\
\hline
\hline
\cite{McGeoch:2013:EEA:2482767.2482797} & \checkmark & & & & & \checkmark & & & & \checkmark \\ % considered tabu instead of SA
\hline
\cite{ibm_blog} & \checkmark & & & & & \checkmark & & & & \\
\hline
\cite{1306.1202} & & \checkmark & & & & \checkmark & & & & \\
\hline 
\cite{PhysRevA.92.042325} & & & \checkmark & & & & \checkmark & \checkmark & & \checkmark \\
\hline 
\cite{king2015performance} & & & \checkmark & & & & \checkmark & \checkmark & & \checkmark \\
\hline
\cite{1604.00319} & \checkmark & & \checkmark & & & & & \checkmark & \checkmark & \checkmark \\
\hline
\cite{PhysRevX.6.031015} & & & & & \checkmark & & \checkmark & & \checkmark & \checkmark \\
\hline
\cite{1701.04579} & & & & \checkmark & & & \checkmark & \checkmark & \checkmark & \checkmark \\
\hline
%OPO study & & & & & & & & & & \\
%\hline
This Work & \checkmark & \checkmark & \checkmark & \checkmark & \checkmark & \checkmark &  & \checkmark & & \checkmark \\
\hline
%%%%%%%%%%%%%%%%%%%%%%%%%%%%%%%%%%%%%
\end{tabular}\\
\vspace{0.1cm}
\caption{A Chronological Summary of IPU Benchmarking Studies}
\vspace{-0.5cm}
\label{tbl:bench_overview}
\end{table*}

\section{A Study of Established Methods}
\label{sec:computations}

This section conducts an in-depth computational study of the established instance generation algorithms and solution methods for IPUs.  The first goal of this study is to understand what classes of problems and parameters are the most challenging, as such cases are preferable for benchmarking.  The second goal is to conduct a validation study of a D-Wave 2X IPU, to clearly quantify its solution quality and runtime performance.  This computational study is divided into two phases. First, a broad parameter sweep of all possible instance generation algorithms is conducted and a commercial mixed-integer programming solver is used to filter out the easy problem classes and parameter settings. Second, after the most challenging problems have been identified, a detailed study is conducted to compare and contrast the three disparate solution methods IP, LNS, and AQC.

Throughout this section, the following notations are used to describe the algorithm results: 
$UB$ denotes the objective value of the best feasible solution produced by the algorithm within the time limit, 
$LB$ denotes the value of the best lower bound produced by the algorithm within the time limit, 
$T$ denotes the algorithm runtime in seconds\footnote{For MIP solvers, the runtime includes the computation of the optimally certificate.}, 
$TO$ denotes that the algorithm hit a time limit of 600 seconds,
$\mu(\cdot)$ denotes the mean of a collection of values,
$sd(\cdot)$ denotes the standard deviation of a collection of values,
 and $max(\cdot)$ denotes the maximum of a collection of values.

\paragraph{Computation Environment:}

% Darwin Doc: cn[300-338,340-371]       2:18:2       128673  baseboard_vendor:HP,cpu_vendor:Intel,cpu_family:broadwell,cpu_model:E5-2695_v4,cpu_base_clock:2.10GHz,numa_nodes:2,clusterondie:no,multithreading:yes,ib:edr,ethernet:1Gb,nvme:no,ssd:no,hdd:no,gpu_count:0

% From Ryan - HPE ProLiant XL170r Gen9 servers with dual Intel E5-2695 v4 2.10GHz CPUs and 128GB of DDR4-2400 memory

The classical computing algorithms are run on HPE ProLiant XL170r servers with dual Intel 2.10GHz CPUs and 128GB memory.  After a preliminary comparison of CPLEX 12.7 \cite{cplex} and Gurobi 7.0 \cite{gurobi}, no significant difference was observed. Thus, Gurobi was selected as the commercial Mixed-Integer Programming (MIP) solver and was configured to use one thread.  The highly specialized and optimized HFS algorithm \cite{HFS_impl_2017} is used as an LNS-based heuristic and also uses one thread.

The IPU computation is conducted on a D-Wave 2X \cite{dwave_2x} adiabatic quantum computer (AQC). This computer has a 12-by-12 chimera cell topology with random omissions; in total, it has 1095 spins and 3061 couplers and an effective temperature of $\bm \tau \in (0.091, 0.053)$ depending on the problem being solved \cite{NIPS2016_6375,1612.05024}.  Unless otherwise noted, the AQC is configured to produce 10,000 samples using a 5-microsecond annealing time per sample and a random gauge transformation every 100 samples.  The best sample is used in the computation of the upper bound value.  The reported runtime of the AQC reflects the amount of time used on the IPU hardware; it does not include the overhead of communication or scheduling of the computation, which adds an overhead of about three seconds.

All of the software used in this benchmarking study is available as open-source via: \ipujson, a language-independent JSON-based Ising model exchange format designed for benchmarking IPU hardware; \ipig, algorithms for IPU instance generation; \ipusolvers, tools for encoding \ipujson~data into various optimization formulations and solvers.\footnote{The source code is available at \url{https://github.com/lanl-ansi/} under the repository names {\sc bqpjson}, {\sc dwig} and {\sc bqpsolvers}.}

\subsection{Identifying Challenging Cases}
\label{sec:computations:cases}

\paragraph{Broad Parameter Sweep:}
In this first experiment, we conduct a parameter sweep of all the inputs to the problem generation algorithms described in Section \ref{sec:tools_gen}.
Table \ref{tbl:all_params} provides a summary of the input parameters for each problem class.  The values of each parameter are encoded with the following triple: $(\mbox{start}..\mbox{stop} : \mbox{step size})$.  When two parameters are required for a given problem class, the cross product of all parameters is used.  For each problem class and each combination of parameter settings, 250 random problems are generated in order to produce a reasonable estimate of the average difficulty of that configuration.  Each problem is generated using all of the decision variables available on the IPU.  The computational results of this parameter sweep are summarized in Table \ref{tbl:grb_all_sweep}.  
%It is important to recall that the FL and FCL generation methods are randomized algorithms and are not guaranteed to succeed.  This explains why the number of cases for these problems is not a multiple of 250.

\begin{table*}[t]
\center
\begin{tabular}{|l||l|l|r||r|r|r|r|r|r|r|r|r|r|r|c|c|}
\hline
Problem & First Param. & Second Param. \\
\hline
\hline
RAN-k & $k \in (1..5:1)$ & \multicolumn{1}{|c|}{NA} \\
\hline
RANF-k & $k \in (1..5:1)$ & \multicolumn{1}{|c|}{NA} \\
\hline
FL-k & $k \in (1..5:1)$ & $\alpha \in (0..1:0.1)$  \\
\hline
FCL-k & $k \in (1..5:1)$ & $\alpha \in (0..1:0.1)$ \\
\hline
WSCN & ${\it wf} \in (-1..1:0.2)$ & ${\it sf} \in (-1..1:0.2)$ \\
\hline
%%%%%%%%%%%%%%%%%%%%%%%%%%%%%%%%%%%%%
\end{tabular}\\
\vspace{0.1cm}
\caption{Parameter Settings of Various Problems.}
%\vspace{-0.5cm}
\label{tbl:all_params}
\end{table*}

The results presented in Table \ref{tbl:grb_all_sweep} indicate that, at this problem size, all variants of the FL, FCL, and WSCN problems are easy for modern MIP solvers.  This is a stark contrast to \cite{PhysRevX.6.031015}, which reported runtimes around 10,000 seconds when applying Simulated Annealing to the WSCN problem.  Furthermore, this result suggests that these problems classes are not ideal candidates for benchmarking IPUs.  In contrast, the RAN and RANF cases consistently hit the runtime limit of the MIP solver, suggesting that these problems are more useful for benchmarking.  This result is consistent with a similar observation in the SAT community, where random SAT problems are known to be especially challenging \cite{Mitchell:1992,Balyo:2017}.  To get a better understanding of these RAN problem classes, we next perform a detailed study of these problems for various values of the parameter $k$. 

\begin{table*}[t]
%\small
%\footnotesize
%\scriptsize
%\tiny
\center
\begin{tabular}{|l||r|r|r||r|r|r|r|r|r|r|r|r|r|r|c|c|}
\hline
Problem & Cases & $\mu(|{\cal N}|)$ & $\mu(|{\cal E}|)$ & $\mu(T)$ & $sd(T)$ & $max(T)$ \\
\hline
% generator: problem_stats.r
% data: sweep-mip_data.csv
% output: cp_paper_all_runtimes.txt
%%%%%%%%%%%%%%%%%%%%%%%%%%%%%%%%%%%%%
%Problem & Cases & mu(bits) & mu(edges) & mu(time) & sd(time) & max(time)
\hline
RAN & 1250 & 1095 & 3061 & \timeout & --- & \timeout \\
\hline
RANF & 1250 & 1095 & 3061 & \timeout & --- & \timeout \\
\hline
FL & 6944 & 1008 & 2126 & 1.82 & 1.06 & 16.80 \\
\hline
FCL & 8347 & 888 & 2282 & 4.19 & 2.81 & 41.40 \\
\hline
WSCN & 30250 & 949 & 2313 & 0.25 & 0.87 & 17.90 \\
\hline
%%%%%%%%%%%%%%%%%%%%%%%%%%%%%%%%%%%%%
\end{tabular}\\
\vspace{0.1cm}
\caption{MIP Runtime on Various IPU Benchmark Problems (seconds)}
%\vspace{-0.5cm}
\label{tbl:grb_all_sweep}
\end{table*}

\paragraph{The RAN and RANF Problems:}
In this second experiment, we focus on the RAN-k and RANF-k problems and conduct a detailed parameter sweep of $k \in (1..10:1)$.  To accurately measure the runtime difficulty of the problem, we also reduce the size of the problem from 1095 variables to 194 variables so that the MIP solver can reliably terminate within a 600 second time limit.  The results of this parameter sweep are summarized in Table \ref{tbl:grb_ran_sweep}.

The results presented in Table \ref{tbl:grb_ran_sweep} indicate that (1) as the value of $k$ increases, both the RAN and RANF problems become easier; and (2) the RANF problem is easier than the RAN problem.  The latter is not surprising because the additional linear coefficients in the RANF problem break many of the symmetries that exist in the RAN problem.  These results suggest that it is sufficient to focus on the RAN-1 and RANF-1 cases for a more detailed study of IPU performance.  This is a serendipitous outcome for IPU benchmarking because restricting the problem coefficients to $\{-1,0,1\}$ reduces artifacts caused by noise and the numeral precision of the analog hardware.

\begin{table*}[t]
%\small
\center
\begin{tabular}{|c|r|r|r||r|r|r||r|r|r|r|r|r|r|r|c|c|}
\hline
$k$ & Cases & $\mu(|{\cal N}|)$ & $\mu(|{\cal E}|)$ & $\mu(T)$ & $sd(T)$ & $max(T)$ & $\mu(T)$ & $sd(T)$ & $max(T)$ \\
\hline
% generator: ran_steps.r
% data: ran-steps-sml-mip_data.csv
% output: cp_paper_ran_k_runtimes.txt, cp_paper_ranf_k_runtimes.txt
%%%%%%%%%%%%%%%%%%%%%%%%%%%%%%%%%%%%%
\hline
\multicolumn{4}{|c||}{Problems of Increasing k} & \multicolumn{3}{c||}{RAN-k} & \multicolumn{3}{c|}{RANF-k}  \\
%%%%%%
%k & Cases & mu(bits) & mu(edges) & mu(time) & sd(time) & max(time)
\hline
1 & 250 & 194 & 528 & 340.0 & 195.0 & \timeout & 14.10 & 15.20 & 82.70 \\
\hline
2 & 250 & 194 & 528 & 89.3 & 64.3 & 481 & 2.97 & 3.41 & 22.70 \\
\hline
3 & 250 & 194 & 528 & 64.8 & 28.3 & 207 & 1.67 & 1.48 & 10.70 \\
\hline
4 & 250 & 194 & 528 & 58.0 & 29.5 & 250 & 1.25 & 0.83 & 6.10 \\
\hline
5 & 250 & 194 & 528 & 49.0 & 23.0 & 131 & 1.12 & 0.77 & 6.98 \\
\hline
6 & 250 & 194 & 528 & 49.0 & 22.4 & 119 & 1.05 & 0.59 & 4.47 \\
\hline
7 & 250 & 194 & 528 & 45.0 & 22.8 & 128 & 1.04 & 0.75 & 7.60 \\
\hline
8 & 250 & 194 & 528 & 44.8 & 23.7 & 121 & 1.01 & 0.62 & 5.43 \\
\hline
9 & 250 & 194 & 528 & 42.3 & 22.3 & 110 & 0.98 & 0.60 & 5.08 \\
\hline
10 & 250 & 194 & 528 & 39.8 & 22.1 & 107 & 0.91 & 0.43 & 3.09 \\
\hline
%%%%%%%%%%%%%%%%%%%%%%%%%%%%%%%%%%%%%
\end{tabular}\\
\vspace{0.2cm}
\caption{MIP Runtime on RAN-k and RANF-k IPU Benchmark Problems (seconds)}
\vspace{-0.8cm}
\label{tbl:grb_ran_sweep}
\end{table*}

\subsection{An IPU Evaluation using RAN-1 and RANF-1}
\label{sec:computations:ipu}

Now that the RAN-1 and RANF-1 problem classes have been identified as the most interesting for IPU benchmarking, we perform two detailed studies on these problems using all three algorithmic approaches (i.e., AQC, LNS, and MIP).  The first study focuses on the scalability trends of these solution methods as the problem size increases, whereas the second study focuses on a runtime analysis of the largest cases that can be evaluated on a D-Wave 2X IPU hardware. 

\paragraph{Scalability Analysis:}
In this experiment, we increase the problem size gradually to understand the scalability profile of each of the solution methods (AQC, LNS, and MIP).  The results are summarized in Table \ref{tbl:algs_comp}.  Focusing on the smaller problems, where the MIP solver provides an optimality proof, we observe that both the AQC and the LNS methods find the optimal solution in all of the sampled test cases, suggesting that both heuristic solution methods are of high quality.  Focusing on the larger problems, we observe that, in just a few seconds, both AQC and LNS find feasible solutions that are of higher quality than what the MIP solver can find in 600 seconds.  This suggests that both methods are producing high-quality solutions at this scale.  As the problem size grows, a slight quality discrepancy emerges favoring LNS over AQC; however, this discrepancy in average solution quality is less than 1\% of the best known value.

\begin{table*}[t]
%\small
\center
\begin{tabular}{|r|r|r||r|r||r|r||r|r|r|r|r|r|r|c|c|}
\hline
 \multicolumn{3}{|c||}{} & \multicolumn{2}{c||}{AQC} & \multicolumn{2}{c||}{LNS} & \multicolumn{3}{c|}{MIP}\\
\hline
Cases & $\mu(|{\cal N}|)$ & $\mu(|{\cal E}|)$ & $\mu(UB)$ & $\mu(T)$ & $\mu(UB)$ & $\mu(T)$ & $\mu(UB)$ & $\mu(LB)$ & $\mu(T)$\\
\hline
% generator: ran_cd.r
% data: ran-cd-mip_data.csv, ran-cd-lns_data.csv, ran-cd-aqc_data.csv
% output: cp_paper_ran_1_runtimes.txt, cp_paper_ranf_1_runtimes.txt
%%%%%%%%%%%%%%%%%%%%%%%%%%%%%%%%%%%%%
\hline
\multicolumn{10}{|c|}{RAN-1 Problems of Increasing Size} \\
%%%%%%
%Cases & mu(bits) & mu(edges) & mu(ub) & mu(time) & mu(ub) & mu(time) & mu(ub) & mu(lb) & mu(time)
\hline
250 & 30 & 70 & -44 & 3.53 & -44 & 10 & -44 & -44 & 0.05 \\
\hline
250 & 69 & 176 & -110 & 3.57 & -110 & 10 & -110 & -110 & 0.48 \\
\hline
250 & 122 & 321 & -199 & 3.60 & -199 & 10 & -199 & -199 & 15.90 \\
\hline
250 & 194 & 528 & -325 & 3.64 & -325 & 10 & -325 & -327 & 340.00 \\
\hline
250 & 275 & 751 & -462 & 3.68 & -462 & 10 & -461 & -483 & \timeout \\
\hline
250 & 375 & 1030 & -633 & 3.73 & -633 & 10 & -629 & -673 & \timeout \\
\hline
250 & 486 & 1337 & -821 & 3.77 & -822 & 10 & -814 & -881 & \timeout \\
\hline
250 & 613 & 1689 & -1038 & 3.77 & -1039 & 10 & -1021 & -1116 & \timeout \\
\hline
250 & 761 & 2114 & -1296 & 3.76 & -1297 & 10 & -1262 & -1401 & \timeout \\
\hline
250 & 923 & 2578 & -1574 & 3.77 & -1576 & 10 & -1525 & -1713 & \timeout \\
\hline
250 & 1095 & 3061 & -1870 & 3.80 & -1873 & 10 & -1806 & -2045 & \timeout \\
\hline
%%%%%%
\hline
\multicolumn{10}{|c|}{RANF-1 Problems of Increasing Size} \\
%%%%%%
%Cases & mu(bits) & mu(edges) & mu(ub) & mu(time) & mu(ub) & mu(time) & mu(ub) & mu(lb) & mu(time)
\hline
250 & 30 & 70 & -53 & 3.53 & -53 & 10 & -53 & -53 & 0.02 \\
\hline
250 & 69 & 176 & -127 & 3.56 & -127 & 10 & -127 & -127 & 0.13 \\
\hline
250 & 122 & 321 & -229 & 3.61 & -229 & 10 & -229 & -229 & 0.67 \\
\hline
250 & 194 & 528 & -370 & 3.66 & -370 & 10 & -370 & -370 & 14.10 \\
\hline
250 & 275 & 751 & -526 & 3.71 & -526 & 10 & -526 & -527 & 128.00 \\
\hline
250 & 375 & 1030 & -719 & 3.76 & -719 & 10 & -719 & -727 & 471.00 \\
\hline
250 & 486 & 1337 & -934 & 3.81 & -934 & 10 & -933 & -954 & 588.00 \\
\hline
250 & 613 & 1689 & -1179 & 3.82 & -1179 & 10 & -1178 & -1211 & \timeout \\
\hline
250 & 761 & 2114 & -1472 & 3.82 & -1472 & 10 & -1470 & -1520 & \timeout \\
\hline
250 & 923 & 2578 & -1786 & 3.82 & -1787 & 10 & -1778 & -1856 & \timeout \\
\hline
250 & 1095 & 3061 & -2121 & 3.86 & -2122 & 10 & -2110 & -2212 & \timeout \\
\hline
%%%%%%%%%%%%%%%%%%%%%%%%%%%%%%%%%%%%%
\end{tabular}\\
\vspace{0.2cm}
\caption{A Comparison of Solution Quality and Runtime as Problem Size Increases on RAN-1 and RANF-1.}
\vspace{-0.8cm}
\label{tbl:algs_comp}
\end{table*}

\begin{figure}[t]
% generator: ran_rts.r
% data: ran-cd12-mip-rts_data.csv, ran-cd12-lns-rts_data.csv, ran-cd12-aqc-rts-0005_data.csv
% output: see pdf file names
\center
    \includegraphics[width=6.2cm]{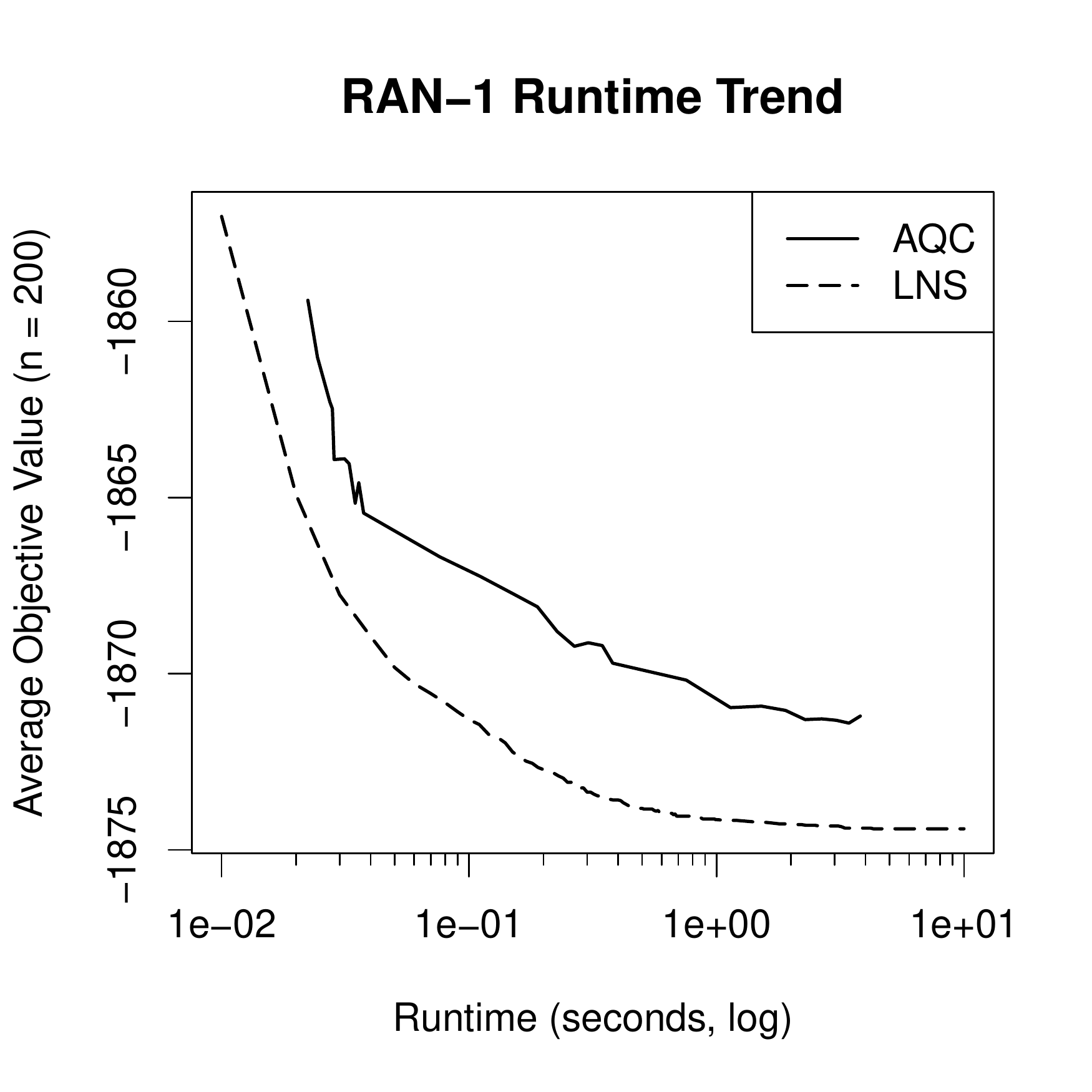} 
    \hspace{-0.5cm}
    \includegraphics[width=6.2cm]{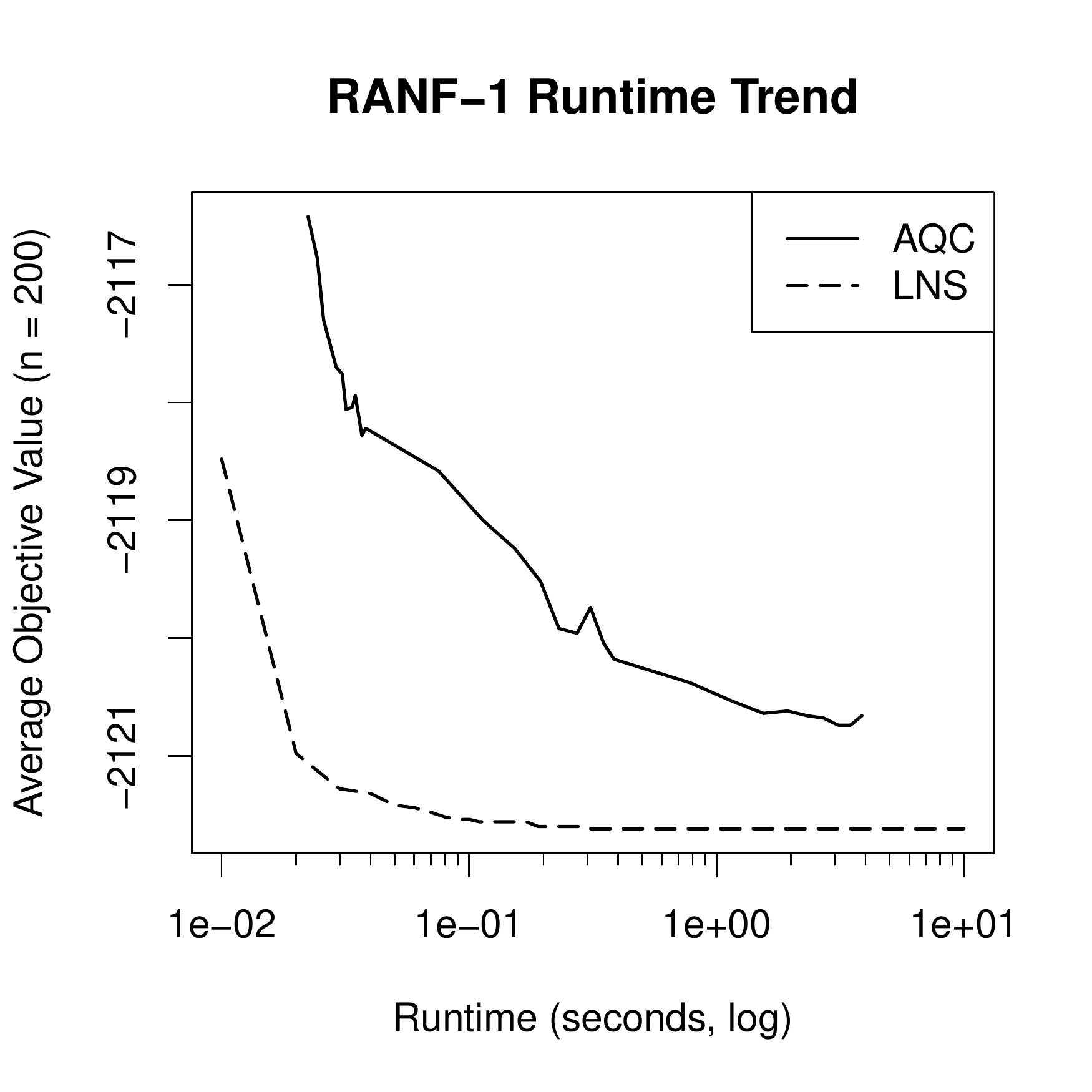} 
\vspace{-0.4cm}
\caption{Detailed Runtime Analysis of the AQC (D-Wave 2X) and LNS Heuristic (HFS) on the RAN-1 (left) and RANF-1 (right) Problem Classes.}
%\vspace{-0.3cm}
\label{fig:algs_runtime}
%\vspace{-0.3cm}
\end{figure}

\paragraph{Detailed Runtime Analysis:}
Given that both the AQC and the LNS solution methods have very similar solution qualities, it is prudent to perform a detailed runtime study to understand the quality vs. runtime tradeoff.  To develop a runtime profile of the LNS algorithm, the solver's runtime limit is set to values ranging from 0.01 to 10.00 seconds.  In the case of the AQC algorithm, the number of requested samples is set to values ranging from 10 to 10,000, which has the effect of scaling the runtime of the IPU process.
%\footnote{It is important to emphasize that the IPU runtime does not include any communication or scheduling overheads.} % cut for space   
The results of this study are summarized in Figure \ref{fig:algs_runtime}.  Note that the stochastic sampling nature of the IPU results in some noise for small numbers of samples.  However, the overall trend is clear.

The results presented in Figure \ref{fig:algs_runtime} further illustrate that (1) the RAN problem class is more challenging than the RANF problem class, and (2) regardless of the runtime configuration used, the LNS heuristic slightly outperforms the AQC; however, the average solution quality is always within 1\% of each other.  Combining all of the results from this section provides a strong validation that even if the D-Wave 2X IPU cannot guarantee a globally optimal solution, it produces high quality solutions reliably across a wide range of inputs. 

% \subsection{Discussion of the RAN Problem}
% [cjc] thinking we should cut this section because it could easily undermine the rest of the computational analysis.
% \begin{enumerate}
%     \item provide an intuition for why this is so heuristics (six regular graph)
%     \item and why this is so hard for complete solvers (i.e. no problem structure / symmetries)
% \end{enumerate}

\section{Conclusion}
\label{sec:conclusion}

This work introduces the idea of Ising processing units (IPUs) as a computational abstraction for emerging physical devices that optimize Ising models.  It highlights a number of unexpected challenges in using such devices and proposes commercial mixed-integer programming solvers as a tool to help improve validation and benchmarking.

A baseline study of the D-Wave 2X IPU suggests that the hardware specific instance generation is a reasonable strategy for benchmarking IPUs.  However, finding a class of challenging randomly generated test cases is non-trivial and an open problem for future work.  The study verified that at least one commercially available IPU is already comparable to current state-of-the-art classical methods on some classes of problems (e.g. RAN and RANF).  
%However, a detailed runtime analysis did not demonstrate any time vs quality configuration where the IPU outperformed a state-of-the-art classical method.  That said, the IPU hardware seems to be largely invariant in quality and runtime across a wide variety of problem classes and instance sizes.  
Consequently, as this IPU's hardware increases in size, one would expect that it could outperform state-of-the-art classical methods because of its parallel computational nature and become a valuable co-processor in hybrid-optimization algorithms.
%This is perhaps not surprising as the IPU hardware is effectively a fully parallelized optimization method.

Overall, we find that the emergence of IPUs is an interesting development for the optimization community and warrants continued study.  Considerable work remains to determine new challenging classes of test cases for validating and benchmarking IPUs.
%and to develop better practices for fairly comparing heterogeneous hardware platforms to established classical computing methods.  
We hope that the technology overview and the validation study conducted in this work will assist the optimization research community in exploring IPU hardware platforms and will accelerate the development of hybrid-algorithms that can effectively leverage these emerging technologies.

%\section*{Acknowledgments}
%We gratefully acknowledge Edward (Denny) Dahl for his feedback on the tools proposed herein and the U.S. Department of Energy for supporting this work through Los Alamos National Laboratory's LDRD Program.

% \begin{enumerate}
%     \item Gurobi
%     \itme 
% \end{enumerate}

%==============================;
%  Include all the references  ;
%==============================;
\clearpage
\bibliographystyle{splncs}
\bibliography{references.bib}

\noindent
%LA-UR-19-22000

\end{document}